\documentclass[12pt]{amsart}
\usepackage[all]{xy}
\usepackage{amssymb}

\newcommand{\cA}{\mathcal{A}}
\newcommand{\cL}{\mathcal{L}}
\newcommand{\cS}{\mathcal{S}}
\newcommand{\bba}{\mathbb{A}}
\newcommand{\bbc}{\mathbb{C}}
\newcommand{\bbf}{\mathbb{F}}
\newcommand{\bbn}{\mathbb{N}}
\newcommand{\bbp}{\mathbb{P}}
\newcommand{\bbq}{\mathbb{Q}}

\newcommand{\bbz}{\mathbb{Z}}
\newcommand{\Pic}{\mathrm{Pic}\,}
\newcommand{\Gal}{\mathrm{Gal}\,}
\newcommand{\ch}{\mathrm{char}\,}
\newcommand{\sep}{\mathrm{sep}}
\newcommand{\prim	}{\mathrm{prim}}
\newcommand{\ord}{\mathrm{ord}}

\newcommand{\fm}{\mathfrak{m}}

\newtheorem{thm}{Theorem}[section]
\newtheorem{cor}[thm]{Corollary}
\newtheorem{lem}[thm]{Lemma}
\newtheorem{prop}[thm]{Proposition}

\newtheorem{defnn}[thm]{Definition}
\newenvironment{definition}{\begin{defnn} \em}{\end{defnn}}
\newtheorem{remarkk}[thm]{Remark}

\newtheorem{examplee}[thm]{Example}

\title{Abelian Varieties over Cyclic Fields}
\author{Bo-Hae Im and Michael Larsen}
\date{\today}
\address{Department of Mathematics, University of Utah, Salt Lake City, Utah
 84112, USA}\email{im@math.utah.edu}
  \address{Department of Mathematics, Indiana University, Bloomington,
Indiana 47405, USA} \email{larsen@math.indiana.edu}

\begin{document}
\begin{abstract} Let $K$ be a field of characteristic $\neq 2$ such that every finite separable
 extension of $K$ is cyclic.
Let $A$ be an abelian variety over $K$.  If $K$ is infinite, then $A(K)$ is Zariski-dense in $A$.
If $K$ is not locally finite, the rank of $A$ over $K$ is infinite.
\end{abstract}
\maketitle

\section{Introduction}
Let $K$ be a field and $K^{\sep}$ a separable closure of $K$.

\begin{definition}
We say $K$ is \emph{cyclic} if $\Gal(K^{\sep}/K)$ is topologically cyclic, i.e. if every
finite separable extension of $K$ is cyclic.
\end{definition}

The standard examples of cyclic fields from our point of view are the quasi-finite fields,
which are perfect fields such that $\Gal(K^{\sep}/K)\cong \hat\bbz$.
Other examples are real-closed fields and separably closed fields.
For any automorphism $\sigma$ of a separably closed field $F$,
the fixed field $F^\sigma$ is cyclic, and all cyclic fields arise in this way.

Our goal in this paper is to understand groups of the form $A(K)$ where $A$ is an abelian variety and $K$ is a cyclic field.
Finite fields are cyclic, and of course the group of points of an abelian variety over a finite field is finite.  Our two main results show that infinite cyclic fields behave quite differently.

\begin{thm}
\label{Density}
Let $A$ be an abelian variety over an infinite cyclic field $K$ with $\ch K\neq 2$.  Then $A(K)$ is Zariski-dense
in $A$.
\end{thm}

Every algebraic extension of a finite field is cyclic, and the above theorem applies, therefore,
when $K$ is an infinite algebraic extension of $\bbf_p$.  In such cases, of course,
$A(K)\subset A(\bar{\bbf}_p)$, which is a torsion group, so the rank of $A$ over $K$ is 0.
When $K$ is \emph{not} algebraic over $\bbf_p$ for any $p$ (i.e., not locally finite), we observe totally different behavior:

\begin{thm}
\label{Rank}
Let $A$ be a non-trivial abelian variety over a cyclic field $K$ which is not locally finite
and such that $\ch K\neq 2$.  Then
$$\dim A(K)\otimes\bbq = \infty.$$
\end{thm}

This theorem has an equivalent formulation as follows:

\begin{thm}
\label{Invariants}
Let $A$ be a non-trivial abelian variety over an infinite, finitely generated field $K_0$ of characteristic $\neq 2$.
Then for all $\sigma\in\Gal(K_0^{\sep}/K_0)$, the group $A((K_0^{\sep})^\sigma)$ has infinite rank.
\end{thm}

In this form, it can be regarded as an unconditional version of an old probabilistic
result of G.~Frey and M.~Jarden \cite{FJ}:

\begin{thm}
Let $A$ be a non-trivial abelian variety over an infinite, finitely generated field $K_0$, and let $n$ be a positive
integer.  The set of $n$-tuples $(\sigma_1,\ldots,\sigma_n)\in \Gal(K_0^{\sep}/K_0)^n$ for which
$A((K_0^{\sep})^{\langle \sigma_1,\ldots,\sigma_n\rangle})$ has infinite rank is of full measure.
\end{thm}

In \cite{La}, the second-named author asked for an unconditional version of this theorem.
The problem remains open for $n>1$.

As a consequence of Theorem~\ref{Rank}, we also can prove the following generalization of Theorem~\ref{Invariants}, which answers in the affirmative Questions 1.1--1.3 in \cite{IL}:

\begin{thm}
\label{Spectrum}
Let $A$ be a non-trivial abelian variety over an infinite, finitely generated
field $K_0$ of characteristic $\neq 2$.  Then for all $\sigma\in\Gal(K_0^{\sep}/K_0)$,
$\sigma$ acts semisimply on $A(K_0^{\sep})\otimes\bbc$.  The eigenvalues are all roots
of unity, and all have countably infinite multiplicity.  
Moreover, the following conditions on $n$ and $\sigma$ are equivalent:
\begin{enumerate}
\item All primitive $n$th roots of unity are eigenvalues.
\item Some primitive $n$th root of unity is an eigenvalue.
\item The closed subgroup $\overline{\langle\sigma^n\rangle}$ is of index $n$ in 
$\overline{\langle\sigma\rangle}$.
\item The group $\overline{\langle\sigma\rangle}$ has an open subgroup of index $n$.

\end{enumerate}

\end{thm}

We remark that prior to this paper, substantial progress was made on the $\dim A = 1$ case
using Heegner point methods.
In particular, the first-named author \cite{im3} settled the case of elliptic curves over $\bbq$.
This was extended to cover many cases over totally real fields and function fields of odd characteristic
\cite{BI}.

The approach of this paper is completely different.
The basic strategy was worked out in the $\dim A=1$
case in \cite{im1}, improving on the ad hoc methods of
\cite{La}.  In rough terms, it goes as follows.
Let $\sigma$ denote a topological generator of $\Gal(K^{\sep}/K)$.
Let $n$ be an even integer and $Y$ an algebraic curve over $K$
which is a Galois $\cA_n$-extension
of $\bbp^1$.  Suppose there exists a $K$-morphism $f\colon Y/ \cA_{n-1}\to A$
whose image generates the abelian variety.
A point $c\in \bbp^1(K)$ which is not in the branch locus of the quotient map
lifts to a $\sigma$-stable $\cA_n$-orbit in $Y(K^{\sep})$; that is,
$\sigma$ acts on the $n$-element set $\{x_1,\ldots,x_n\}$ of preimages of $c$
in $Y/\cA_{n-1}$ by
an even permutation.  In particular, the $\sigma$-orbit
$$\{x_1,\sigma(x_1),\sigma^2(x_1),\ldots,\sigma^{m-1}(x_1)\}$$
is a \emph{proper} subset of $\{x_1,\ldots,x_n\}$, since an $n$-cycle is an odd permutation.
The sum
$$f(x_1)+f(\sigma(x_1))+\cdots+f(\sigma^{m-1}(x_1))$$
gives a point in $A(K)$.
We show that points of this form are Zariski-dense
whenever $K$ is infinite and form a group of infinite rank whenever $K$ contains a Hilbertian subfield over which $A$ can be defined.

The above sketch is in fact an oversimplification.
In reality, we do not know how to find the desired curve $Y$ and morphism $f$ over $K$; instead, we make a separable extension of scalars in order to find the requisite geometric data and then use
what amounts to a restriction of scalars argument to control the traces of the resulting points.

Even over a separably closed field it is not trivial to find $Y$.  The idea is to start with a non-singular
projective curve $X\subset A$ and consider separable degree $n$ morphisms $\pi\colon X\to \bbp^1$
such that the Galois group of the splitting field $M$ of the function field of $X$ over the function field
of $\bbp^1$ is $\cS_n$.  If $Y$ denotes the projective non-singular curve with function field $M$, then
$X = Y/\cS_{n-1}$.  By carefully controlling the ramification of $\pi$, we can arrange that
$Y/\cA_n$ has genus $0$.  Then $Y/\cA_{n-1}$ is a double cover of $X$ and therefore maps to $A$
with image $X$.

This paper contains three substantial improvements over \cite{im1}.   The first is to work with general abelian varieties rather than elliptic curves.  This means, in effect, that we must find rational curves on
varieties of the form $A^n/\cA_n$.  We do it by choosing $X\subset A$ and using Riemann-Roch arguments to find rational curves on $X^n/\cA_n$.  The second is to treat the case of finite (odd) characteristic together with the characteristic zero case.  This requires an analysis of wild ramification, which in turn depends on some harder group theory.  The third is the introduction of a trace argument, which frees us from the difficulty of finding $K$-rational data.  (This difficulty led to extra hypotheses and additional work in \cite{La}, \cite{im1}, and \cite{im4}.)  This third improvement was facilitated by the shift in viewpoint from fields of the form
$(K_0^{\sep})^\sigma$ (where $K_0$ is Hilbertian) to abstract cyclic fields $K$---at the one point where Hilbert irreducibility is used (and it is now used in a very weak form), we pass to a finitely generated subfield.

\section{Constructing $\cA_n$-extensions}

Throughout this section, we work over a separably closed field $F$ of characteristic
$p\ge 0$, $p\neq 2$.

Our goal in this section is to show that every abelian variety $A$ belongs to a diagram
\begin{equation}
\label{big-diagram}
\xymatrix{Y\ar[d] \\ Y/\cA_{n-1}\ar[d]_-\phi\ar[r]^-f&A\\ Y/\cA_n\ar@{=}[r]& \bbp^1\\}
\end{equation}
where $Y$ is a non-singular projective curve and $f$ is a non-constant map.

Let $A$ be a non-trivial abelian variety over $F$.  By Bertini's theorem, we can choose
a non-singular projective curve $X$ on $A$.  Translating if necessary, we may assume
the identity $0$ lies on $X$.
Let $g$ denote the genus of $X$.

We fix $k > 2g$ and choose a point
$x\in X(F)$, $x\neq 0$.   If $g=1$, we insist that $[x]-[0]$ is not a 2-torsion class in $\Pic X$; otherwise, there exist finitely many non-trivial involutions of $X$ and we assume that
$x$ is not fixed by any of them.
Let $P=p$ and $Q = \frac{p+1}2$ if $p>0$ and otherwise $P=Q=1$.
By Riemann-Roch, there exists an effective divisor
$D$ of degree $k$ such that
$$2Q D+[x]-(2Q k-2g+3)[0]$$
is in the canonical divisor class.
We assume that
$$D=[w_1]+[w_2]+\cdots+[w_k]$$
is a sum of $k$ distinct points.
This can always be arranged:

\begin{lem}
Every divisor class of degree $k  > 2g$ in $X$ defined over $F$
contains an effective divisor which is the
sum of $k$ distinct points in $X(F)$.
\end{lem}

\begin{proof}
By Riemann-Roch, the fibers of the summation map $X^k\to \Pic X$ are
all of dimension $k-g$.  It therefore suffices to show that the fibers of the
map $f\colon X^{k-1}\to \Pic X$ defined by
$$f(x_1,\ldots,x_{k-1}) = 2x_1+x_2+\cdots+x_{k-1}$$
are all of dimension $k-g-1$.  This is so because for each fixed value of $x_1$,
the summation map on the remaining $k-2$ variables has all fibers isomorphic to
a projective space of dimension $k-2-g$.
\end{proof}

Let $\omega$ denote a meromophic differential form
on $X$, whose divisor is
$$2Q D + [x] - (2Q k-2g+3)[0].$$
We consider the map $j\mapsto j^{2P}\omega$ sending
$$H^0(X,\cL(m[0]))\to H^0(X,\cL((2P m+2Q k-2g+3)[0] - [x] - 2Q D)).$$
We also consider the exterior derivative
\begin{multline*}
d\colon H^0(X,\cL((2P m+2Q k-2g+2)[0])) \\
\to H^0(X,\Omega^1_X\otimes\cL((2P m+2Q k-2g+3)[0])).
\end{multline*}
Its kernel is $1$-dimensional, so by Riemann-Roch, its image is of codimension $\le 2g$.
If $F=\bbc$, the image consists of the meromorphic differential forms, holomorphic except for a pole
of order at most $2P m+2Q k-2g+3$ at $0$, for which the period integrals are all zero.
The condition on $j$ that
\begin{equation}
\label{integ}
dh = j^{2P}\omega
\end{equation}
has a solution $h\in H^0(X,\cL(m[0]))$
defines a system of $\le 2g$ degree-$2P$ forms.   The dimension of the associated projective
variety $Z_m$ is therefore at least $m-3g$.  Moreover, $Z_{m+1}$ is obtained by intersecting
$Z_m$ by the codimension $\le 1$ projective subspace consisting of divisors which
have a pole at $0$ of order $\le m$.  Therefore, $\dim Z_{m+1} \le \dim Z_m+1$.
It follows that equality holds for all $m\gg 0$.  In fact,
given $k$, there are at most $3g$ values of $m$
for which we cannot choose $j$ to have a pole of exact order $m$.

A function $j$ in
$$H^0(X,\cL(m[0]))\setminus H^0(X,\cL((m-1)[0]))$$
for which (\ref{integ}) has a solution
determines a morphism $h\colon X\to\bbp^1$ of degree
$$2q := 2P m+2Q k-2g+2$$
up to additive translation.
We would like $m$ to be small compared to $k$.  At the same time,
we would like $q$ to be a sufficiently large prime.
The following lemma shows that both goals can be achieved
at once:

\begin{lem}
\label{choose}
Let $g$ be a positive integer, $\epsilon > 0$, and $N>0$.
Given fixed integers $P$ and $Q$ which are relatively prime and a set
$S\subset \bbn^2$ such that for each $k$,
$$|\{m\in \bbn\mid (m,k)\in S\}| \le 3g,$$
we can find a pair $(m,k)\notin S$ such that
$$q = P m + Q k - g + 1 > N$$
is prime and $m < q^\epsilon$.
\end{lem}

 \begin{proof}
 For each $M$, the number of pairs $(k,m)$ such that $P m + Q k - g + 1$
 is a prime between $M/2$ and $M$ and $m < (M/2)^\epsilon$
 grows like $M^{1+\epsilon}/\log M$.  On the other hand, the subset of these
 pairs belonging to $S$ is $O(M)$.
  \end{proof}

Thus, for some sufficiently large prime $q$, we can construct a morphism
$h \colon X\to\bbp^1$ of degree $2q$.
We can normalize so that $h (x) = 0$ and $h (0)=\infty$.
The resulting morphism has ramification degree $2q$ at $0$.  It has ramification
degree
$$2+ 2P\,\ord_x j\le 2+2P m$$
at $x$, degree
$$1+2Q+2P\,\ord_{w_i} j$$
for each point $w_i$ appearing in the divisor $D$, and degree
$$1+2P\,\ord_z j$$
for all other points $z$.  In no case does $p$ divide the ramification degree, and except at $x$ and at
$0$, the ramification degree is always odd.  By Lemma~\ref{choose}, we may further assume that the ramification degree at $x$ is less than $\frac{\sqrt{2q}-1}{2}$.  Setting $x_1=x$, all the hypotheses
of the following proposition are therefore satisfied.

\begin{prop}
\label{largeq}
Let $q$ be an odd prime not equal to the characteristic $p\ge 0$ of $F$ and larger than some
sufficiently large absolute constant.
Let $X$ be a projective non-singular curve and $h \colon X\to \bbp^1$ a morphism
of degree $2q$.   For each point $c\in \bbp^1(F)$, let $D_c$ denote $h^{-1}(c)$,
regarded as a divisor on $X$.  Then the splitting field of the function field of $X$ over
that of $\bbp^1$ has Galois
group $\cS_{2q}$ if all of the following conditions hold:

\renewcommand{\theenumi}{\roman{enumi}}
\begin{enumerate}
\item $D_\infty = 2q[0]$ for some point $0\in X(F)$.
\item $D_0$ is of the form $\sum a_i[x_i]$ where $a_1$ is even, $2a_1 < \sqrt{2q}-1$,
and the other coefficients $a_i$ are prime to $2P$.
\item For all $c\notin\{0,\infty\}$,
$D_c=\sum b_{c,i}[y_{c,i}]$ where all $b_{c,i}$ are prime to $2P$.
\item There is no non-trivial involution of $X$ fixing $0$ and $x_1$.
\end{enumerate}
\renewcommand{\theenumi}{\arabic{enumi}}

\end{prop}

\begin{proof}
Let $K_{\bbp^1}$ and $L_X$ denote the function fields of $\bbp^1$ and $X$ respectively, and let
$K_{\bbp^1}^{\sep}$ denote a separable closure of the former.
As $\deg h $ is not divisible by $p$,
it follows that $L_X$ is separable over $K_{\bbp^1}$.  Let $M$ denote the splitting field
of $L_X$ in $K_{\bbp^1}^{\sep}$.
Thus $\Gal(M/K_{\bbp^1})$ acts faithfully on the $2q$-element set $E$
of $K_{\bbp^1}$-embeddings of $L_X$ in $K_{\bbp^1}^{\sep}$.

For each point $c\in \bbp^1(F)$, we let $K_c$ denote the completion of $K_{\bbp^1}$ at $c$,
and we fix a separable closure $K_c^{\sep}$ and an embedding
$\iota\colon K_{\bbp^1}^{\sep}\hookrightarrow K_c^{\sep}$ which extends
$K_{\bbp^1}\hookrightarrow K_c$.
If $D_c = \sum_i b_{c,i}[y_{c,i}]$, then the usual prime decomposition theory for Dedekind domains  gives
$$L_X \otimes_{K_{\bbp^1}} K_c = \prod_i L_{y_{c,i}},\ [L_{y_{c,i}}:K_c] = b_{c,i}.$$
Every $K_{\bbp^1}$-embedding of $L_X$ in $K_c^{\sep}$ factors through a unique
$L_{y_{c,i}}$, and every element of $\Gal(K_c^{\sep}/K_c)$ preserves this factor.
This defines a natural equivalence relation $\sim_c$ on the $2q$-element set of such embeddings,
and the choice of inclusion $\iota$ identifies this set with the set $E$.  The equivalence classes
have cardinality $b_{c,i}$.  Moreover, $\Gal(K_c^{\sep}/K_c)$ acts transitively on the $b_{c,i}$-element subset of $E$ consisting of embeddings which factor through $L_{y_{c,i}}$
since it acts transitively on the set of $K_c$-embeddings of $L_{y_{c,i}}$ in $K_c^{\sep}$.

For any point $c$, we define $I_c$ to be the decomposition group of $c$ in
$\Gal(M/K_{\bbp^1})\subset \cS_{2q}$, i.e. the image under the composition
of $\iota$ and $\Gal(K_{\bbp^1}^{\sep}/K_{\bbp^1})\to \Gal(M/K_{\bbp^1})$ of
$\Gal(K_c^{\sep}/K_c)$.
Note that as $F$ is separably closed,
there is no distinction between decomposition group and inertia group.  If $W_c$ denotes the wild ramification group of $c$, then $W_c$ is a $p$-group (i.e., trivial if $p=0$) which is normal
in $I_c$, and such that $I_c/W_c$ is cyclic.   We have seen that $I_c$ acts transitively
on each $b_{c,i}$-element $\sim_c$-equivalence class of $E$.

Every orbit of a $p$-group is of $p$-power order, so if $G$ is a group acting transitively
on a set of prime-to-$p$ order, a $p$-group normal in $G$ must be trivial.  By hypotheses
(i)--(iii), no $b_{c,i}$ is divisible by $p$, so
$W_c$ is trivial for all $c\in\bbp^1(F)$.   Thus, $I_c$ is topologically cyclic and acts on $E$
by a permutation consisting of cycles of lengths $b_{c,i}$.

Next we assert that the action of $\Gal(M/K_{\bbp^1})$ on $E$ is primitive.  If not,
it must respect a partition into two $q$-element subsets or a partition into $q$ pairs, corresponding
to the cases that $h$ factors as a composition $X\to X'\to \bbp^1$ where
$X'\to \bbp^1$ has degree $2$ or $q$ respectively.
The first is impossible because elements of type (ii) and (iii) can respect such a partition only
by fixing each $q$-element subset, in which case $X'\to\bbp^1$ would be a double cover ramified only over $\infty$, which is impossible in characteristic $p\neq 2$.

For the second, let $E'$ denote the set of $K_{\bbp^1}$-embeddings of $L_{X'}$ in
$K_{\bbp^1}^{\sep}$, so each of the $q$ elements of $E'$ extends to a pair of elements of
$E$.  If $\sigma$ is a permutation of $\{1,2,\ldots,2q\}$ which fixes the set of pairs
$\{\{1,q+1\},\{2,q+2\},\ldots,\{q,2q\}\}$, $\sigma'$ is the image of $\sigma$ in $\cS_q$, and
$i,\sigma'(i),{\sigma'}^2(i),\ldots,{\sigma'}^{k-1}(i)$ is an orbit in $\sigma'$, then $\sigma$ determines a permutation on
the set
\begin{equation}
\label{pre-orbit}
\{i,\sigma'(i),{\sigma'}^2(i),\ldots,{\sigma'}^{k-1}(i),q+i,q+{\sigma'}(i),\ldots,q+{\sigma'}^{k-1}(i)\}
\end{equation}
which is either a single $2k$-cycle or a pair of $k$-cycles.  If $\sigma$ is a generator of inertia
at $c$ acting on $E$, then $\sigma'$ is the same generator acting on $E'$; a $\sigma'$-orbit corresponds
to a point $x'\in X'$ lying over $c$.
The restriction of $\sigma$ to  (\ref{pre-orbit}) consists of one cycle for each $x$ lying over
$x'$, so it consists of two $k$-cycles if and only if $X\to X'$ is unramified at $x'$.
It follows that an $E$-orbit whose order is different from that of any other $E$-orbit
must be associated to the case that $X\to X'$ is ramified at $x'$.  Applying this with hypotheses
(i) and (ii), we see that $0$ and $x'$ both lie over ramification points of the double cover
$X\to X'$ and they are therefore fixed point of the corresponding involution, contrary to
hypothesis (iv).

By (ii), $\Gal(M/K_{\bbp^1})$ contains a non-identity element, a suitable odd power of a generator of $I_0$, which fixes a subset of $E$
whose complement has order $\le a_1\le (\sqrt{2q}-1)/2$.
By a theorem of Babai \cite{Ba}, for $q\gg 0$, this implies
$\cA_{2q}\subset \Gal(M/K_{\bbp^1})$.  By (ii), $\Gal(M/K_{\bbp^1})=\cS_{2q}.$

\end{proof}

\begin{cor}
Let $Y/F$ denote the projective non-singular curve with function field $M$ (the splitting field of $L_X$
over $K_{\bbp^1}$.)
Then $\cA_n$ acts on $Y$, and there exists a diagram of type (\ref{diagram}).

\end{cor}

\begin{proof}
Let $f$ denote the (non-constant) composition of
$$Y/\cA_{n-1}\to Y/\cS_{n-1} = X$$
and $X\hookrightarrow A$, and let $\phi$ denote the quotient map
$Y/\cA_{n-1}\to Y/\cA_n$.  We have the following diagram:
\begin{equation}
\label{diagram}
\xymatrix{Y\ar[d] & & \\
Y/\cA_{n-1}\ar[d]_-\phi\ar[r]\ar@/^1.5pc/[rr]^f &X\ar[r]\ar[d]_h&A \\
Y/\cA_n\ar[r]&Y/\cS_n\ar@{=}[r]& \bbp^1\\ }
\end{equation}
All that remains is to check that $Y/\cA_n$ is of genus $0$.
However, $Y/\cA_n$ is a double cover of $Y/\cS_n = \bbp^1$ which can
only be ramified over $0$, $\infty$, and the values of $h$ at the points of $D$
and the zeroes of $j$.  From condition (iii), however, we know that all elements of
the inertia group of $c\notin\{0,\infty\}$ act as even permutations of $E$.
It follows that this double cover is ramified only at $0$ and $\infty$, and as $p\neq 2$,
this implies $Y/\cA_n$ is of genus $0$.
\end{proof}

\section{The Density Theorem}

\begin{prop}
\label{Non-constant}
Let $F$ be a separably closed field.
For every integer $n\ge 3$, every diagram (\ref{diagram}) defined over $F$, every integer $0<m<n$,
and every finite subset $S\subset A(F)$,
there exists a non-empty open subset $U\subset Y/\cA_n(F)$ such that for $c\in U$
and any $m$ distinct points $x_1,\ldots,x_m\in Y/\cA_{n-1}(F)$ mapping to
$c$, the sum $f(x_1)+\cdots+f(x_m)\notin S$.
\end{prop}

\begin{proof}
Let $C$ denote the set of left cosets of $\cA_{n-1}$ in $\cA_n$.
For $c\in Y/\cA_n$, the divisor $D_c$ which gives the preimage of $c$ in
$Y/\cA_{n-1}$ has divisor class independent of $c$.  Therefore, the sum of $f(x_i)$
over all $x_i$ in $D_c$ (taken with multiplicities) is constant in $A$.
It follows that we can deduce the case $m=n-1$ from the complementary case
$m=1$.  We therefore assume that $m\le n-2$, which means that
$\cA_n$ acts transitively on all ordered $m$-element subsets of $C$.

Every element of $C$ determines a morphism $Y\to Y/\cA_{n-1}$.
An ordered $m$-tuple of distinct $\cA_{n-1}$-cosets
of $L_{Y/\cA_{n-1}}$ in $K_X^{\sep}$ determines a morphism
$i_m\colon Y\to (Y/\cA_{n-1})^m$.  By transitivity, any $m$ distinct elements of $Y/\cA_{n-1}$
with the same image in $Y/\cA_n$ lie in the image of $Y$.  If the map
$Y\to A$ obtained by composing $i_m$ with
$$(x_1,\ldots,x_m)\mapsto f(x_1)+\cdots+f(x_m)$$
is constant, then for all $c_i$ in the complement $V\subset Y/\cA_n$ of the ramification
locus of $Y\to Y/\cA_n$ and all
$x_{i1},x_{i2}\in D_{c_i}$, we have $f(x_{i1})=f(x_{i2})$.   As
$$\sum_{j=1}^n f(x_{ij}) = n f(x_{i1})$$
is independent of $c_i$, it follows that up to $n$-torsion, $f(x)$ is constant on the preimage of
$V$ in $Y/\cA_{n-1}$ and therefore constant on the whole curve, contrary to assumption.

\end{proof}

\begin{lem}
\label{elementary-density}
Let $E$ be any infinite field and $F\subset E^{\sep}$ a finite extension of $E$.
Let $\sigma_i\colon F\to E^{\sep}$, $i=1,\ldots,k$, denote a finite sequence of mutually distinct
$E$-embeddings of
$F$ in $E^{\sep}$.  Then the set
\begin{equation}
\label{dense-set}
\{(\sigma_1(x),\ldots,\sigma_k(x))\mid x\in F\}
\end{equation}
is Zariski-dense in $\bba^k$.
\end{lem}

\begin{proof}
It suffices to prove the lemma in the case that the set $\{\sigma_1,\ldots,\sigma_k\}$
consists of all $E$-embeddings of $F$ in $E^{\sep}$.
As the $\sigma_i$ are linearly independent, there exists
an invertible linear transformation of $(E^{\sep})^k$, i.e. an $E^{\sep}$-automorphism
of $\bba^k$, which maps the set (\ref{dense-set}) to $F^k$; as $F$ is infinite,
$F^k$ is dense in $\bba^k$.

\end{proof}

We can now prove Theorem~\ref{Density}.

\begin{proof}
Without loss of generality, we may and do assume that $A$ is simple over $K$.
If $A(K)$ is not dense in $A$, it must therefore be finite.
We fix a prime $q$ large enough that Proposition~\ref{largeq} applies, and let $n=2q$.  Thus
we obtain a diagram of type (\ref{diagram}) defined over $K^{\sep}$.  Let $L\subset K^{\sep}$
denote a finitely generated $K$-subfield over which (\ref{diagram}) can be defined.
As $K$ is cyclic, $L$ is a finite cyclic extension.  Let $k=[L:K]$.

Let $c$ be an element of $L\subset \bbp^1(L)=Y/\cA_n(L)$.  We assume $c$
is not a branch point of
$g\colon Y/\cA_{n-1}\to Y/\cA_n$
and choose $x_1\in \phi^{-1}(c)$.
The $\sigma$-orbit of $c$ has length $d$, a divisor of $k$, and the $\sigma$ orbit of $x_1$ has
length $dm$, where $m$ is the length of the $\sigma^d$-orbit of $x_1$, i.e., the cardinality of the intersection of the $\sigma$-orbit of $x_1$ with the pre-image $D_c$ of $c$ in $Y/\cA_{n-1}(K^{\sep})$.
As $\sigma^d$ fixes $c$, it induces an even permutation of this intersection, so it cannot be all of
$D_c$.

Clearly
\begin{equation}
\label{big-trace}
\sum_{i=0}^{md-1} f(\sigma^i(x_1)) = \sum_{i=0}^{d-1}\sum_{j=0}^{m-1} f(\sigma^{jd}(\sigma^i(x_1)))
\end{equation}
is a $\sigma$-invariant element of $A(K^{\sep})$ and therefore is an element of $A(K)$.
We claim that for every finite subgroup $G\subset A(K)$,
there exists $x_1$ such that the expression in
(\ref{big-trace}) is not
in $G$.  To prove this we use the right hand side of (\ref{big-trace}).  This is a sum of the form
\begin{equation}
\label{double-sum} \sum_{i=0}^{d-1} \sum_{j=0}^{m-1} f(x_{i,j})
\end{equation}
where $x_{i,0},x_{i,1}\ldots,x_{i,m-1}$ are distinct elements of the
fiber $D_{c_i}$ for some $c_i\in  K^{\sep}$, $i=0,\ldots,d-1$.  As
the $d$-tuples $(c,\sigma(c),\ldots,\sigma^{d-1}(c))$ are dense in
$\bba^d$, it suffices to prove that there exists \emph{some}
$d$-tuple $(c_0,\ldots,c_{d-1})$ such that for all choices $x_{i,j}$
as above, the expression (\ref{double-sum}) is not in $G$.  This is
an immediate consequence of Proposition~\ref{Non-constant}.

\end{proof}

\section{The Infinite Rank Theorem}
In this section, we assume that $K$ is not algebraic over any finite field.  This implies that $K$ contains
a subfield $K_0$ with the same property
which is finitely generated over a prime field and such that $A$ is obtained by extension of scalars
from an abelian variety $A_{K_0}$ defined over $K_0$.  Such a field $K_0$ is Hilbertian
\cite[Chap.~9,~Th.~4.2]{Lang}.  We use only the following, weak consequence of this fact, which follows from \cite[Chap.~9,~Prop.~3.3]{Lang}:
For every finite extension $L_0$ of $K_0$, there exists a quadratic extension of $K_0$ which is not contained in $L_0$.

We will need the following proposition, which was proved by J.~Silverman \cite{sil}
in the case of elliptic curves over number fields.  Generalizing to higher dimensional abelian varieties
is trivial, but generalizing to finitely generated fields in arbitrary characteristic requires extra work.

\begin{prop}
\label{torsion-bound}
Let $K_0$ be a field finitely generated over a prime field, $A_{K_0}$ an abelian variety over
$K_0$, and $n$ an integer.
Then the set of all points of $A_{K_0}$ which are
defined over extensions of $K_0$ of degree $\le n$ is finite.
\end{prop}

\begin{proof}
Any such degree $\le n$ extension of $K_0$ is contained in a separable degree $\le n$
extension of some purely inseparable extension $K'_0$ of degree $\le n$ over $K_0$.
As $K_0$ is finitely generated, there exists a finitely extension
$K''_0$ containing all such $K'_0$; replacing $K_0$ with $K''_0$, without loss of generality
we may consider only separable extensions of $K_0$.

As $A_{K_0}$ is defined over $K_0$, it is defined over some integral domain $R$ finitely generated
over $\bbz$ whose field of fractions is $K_0$; inverting some non-zero element of $R$ if necessary,
we may assume without loss of generality that $A_{K_0}$ extends to an abelian scheme over $R$.
Every finitely generated $\bbz$-algebra has a maximal ideal $\fm$ with finite residue field $\bbf_r$.
If $\ell$ denotes a rational prime distinct from the characteristic of $\bbf_r$, then
there exists a Frobenius element $\mathrm{Frob}_{\fm}\in \Gal(K_0^{\sep}/K_0)$ whose eigenvalues
$\lambda_1,\ldots,\lambda_{2g}$ on $H^1(A_{K_0^{\sep}},\bbz_\ell)$
satisfy $|\iota(\lambda_i)| = r^{1/2}$
for all embeddings $\iota\colon \bbq_\ell\to \bbc$.  Thus, $\prod_i (\lambda_i^{n!} - 1)$
is a non-zero rational integer, bounded by $(r^{n!/2}+1)^{2g}$.  If $\ell$ does not divide this integer,
then $\mathrm{Frob}_{\fm}^i$ cannot fix any non-zero element of $H^1(A_{K_0^{\sep}},\bbf_\ell)$ for $1\le i\le n$.
It follows that there is no rational $\ell$-torsion point over any extension of $K_0$ of degree $\le n$.

Thus, it suffices to bound $\ell$-power torsion over degree $\le n$ extensions of $K_0$
for finitely many primes $\ell$, including the characteristic of $\bbf_r$.  Fix one.
Choose $M$ such that $\ell^M > n!$.
Replacing $K_0$ by a finite extension, if necessary, and applying the Mordell-Weil theorem for
finitely generated fields \cite{Ner}, we may assume that for some $N>M$
\begin{equation}
\label{l-power}
 A_{K_0}(K_0^{\sep})[\ell^M]\subset A_{K_0}(K_0)\cap A_{K_0}(K_0^{\sep})[\ell^\infty]\subset A_{K_0}(K_0^{\sep})[\ell^N].
\end{equation}
Let $L_0$ denote the splitting field of any separable degree $\le n$ extension of $K_0$
and $G=\Gal(L_0/K_0)$.
Suppose $P$ is a point of $A_{K_0}$ of order $\ell^{M+N}$ which is defined over  $L_0/K_0$.
If $\ell^{M-1}(\sigma(P)-P) = 0$ for all $\sigma\in G$ then
$\ell^{M-1}P\in A_{K_0}(K_0)$ is an $\ell$-torsion point of order $\ell^{N+1}$, contrary to
(\ref{l-power}).
Multiplying $P$ by $\ell^s$ for some $s\ge 0$,
we may assume that  $\sigma(P)-P\in A_{K_0}(L_0)[\ell^M]$ for all $\sigma$
and  $\sigma(P)-P\notin A_{K_0}(L_0)[\ell^{M-1}]$ for some $\sigma$.  Thus
$\sigma\mapsto \sigma(P)-P$ is a homomorphism $G\to A_{K_0}(L_0)[\ell^M]$
 of order $\ell^M$.  This is
 impossible since $|G| \le n! < \ell^M$.  We conclude that
$A_{K_0}(L_0)\subset  A_{K_0}[\ell^{M+N-1}]$.
 \end{proof}

We now prove Theorem~\ref{Rank}.

\begin{proof}
Let $A$ be an abelian variety over $K$.   Let $K_0$ be a Hilbertian subfield of $K$
such that $A$ is obtained by extension of scalars from an abelian variety $A_{K_0}$
over $K$.  We fix a prime $q$ large enough that Proposition~\ref{largeq} applies, set $n=2q$,
and let $L_0$ be a finite Galois extension of $K_0$ such that a diagram of type (\ref{diagram}) can be defined over $L_0$.

We assume given a finite set $S=\{a_1,\ldots,a_s\}$ of points of $A_{K_0}$ defined
over subfields $K_1\subset\ldots\subset K_s$ of $K$ which are finite extensions of $K_0$.
Our goal is to find a point $a_{s+1}\in A_{K_0}(K\cap K_s^{\sep})$
and an element $\tau\in \Gal(K_s^{\sep}/K_s)$ such that $\tau(a_{s+1}) - a_{s+1}$ is non-torsion.
Then any $\bbz$-linear relation
$$n_1 a_1 + \cdots + n_s a_s + n_{s+1} a_{s+1} = 0$$
implies
$$n_1 a_1 + \cdots + n_s a_s + n_{s+1} \tau(a_{s+1}) = 0$$
and therefore $n_{s+1}=0$, so $a_{s+1}$ is independent of the previous $a_i$
in $A(K)\otimes\bbq$.   Our goal, then, is to find a finite separable $K_s$-extension
$K_{s+1}\subset K$, a point $a_{s+1}\in A_{K_0}(K_{s+1})$,
and $\tau\in \Gal(K_s^{\sep}/K_s)$ such that $\tau(a_{s+1})-a_{s+1}$ is not torsion.
Note that if this difference \emph{is} torsion, then by Proposition~\ref{torsion-bound}, it lies in a finite
subgroup $T_{s+1}$ of $A_{K_0}(K_s^{\sep})$ which does not depend on the choice of $a_{s+1}$.

Let $L_s = K_s L_0$, which is a finite Galois extension of $K_s$.
Let $K_s(\sqrt{D})$ denote a quadratic extension of $K_s$ whose intersection
with $L_s$ is $K_s$.   We have already remarked that this is possible by Hilbert irreducibility.
Choose $\tau\in\Gal(K_s^{\sep}/K_s)$ such that $\tau(\sqrt D) = -\sqrt D$ and $\tau$ acts trivially
on $L_s$.
Let $M_s = L_s(\sqrt D)$, and choose $c\in M_s\subset  Y/\cA_n(M_s)$,
which we assume is not a branch point of $g$.   Choose $x_1\in \phi^{-1}(c)$.

Let $d$ denote the order of the image of $\sigma$ in $\Gal(M_s/K_s)$.
As in the proof of Theorem~\ref{Density},
the $\sigma^d$-orbit of $x_1$ has $m<n$ points.  We claim that $x_1$ can be chosen so that
$$\sum_{i=0}^{md-1}f(\sigma^i(x_1)) - \sum_{i=0}^{md-1}f(\tau\sigma^i(x_1)) \notin T_{s+1}.$$

By Lemma~\ref{elementary-density}, the set of ordered $2d$-tuples
$$\{(c,\sigma(c),\ldots,\sigma^{d-1}(c),-\tau(c),-\tau\sigma(c),\ldots,-\tau\sigma^{d-1}(c))\mid
c\in M_s\}$$
is Zariski-dense in $\bba^{2d}$.
It therefore suffices to show that there exists a $2d$-tuple
$(c_0,c_1,\ldots,c_{2d-1})\in(K^{\sep})^{2d}$ such that for all $x_{i,j}\in Y/\cA_{n-1}(K^{\sep})$
for which $x_{i,0},\ldots,x_{i,m-1}$ are pairwise distinct elements of the fiber $\phi^{-1}(c_i)$
for $i=0,1,\ldots,2d-1$, we have
$$\sum_{i=0}^{2d-1}\sum_{j=0}^{m-1} f(x_{i,j}) \notin T_{s+1}.$$
This follows from Proposition~\ref{Non-constant}.
\end{proof}

\section{Spectra of Automorphisms}

In this section, we prove Theorem~\ref{Spectrum}.  
We begin by proving
that $\sigma$ is always diagonalizable and its eigenvalues are roots of unity:

\begin{lem}
If $K_0$ is a finitely generated field, $A/K_0$ is an abelian variety, and $\sigma\in \Gal(K_0^{\sep}/K_0)$, then the action of
$\sigma$ on $A(K_0^{\sep})\otimes\bbc$ is semisimple, and all the eigenvalues are roots of unity.
\end{lem}

\begin{proof}
Every vector $v$ in $A(K_0^{\sep})\otimes\bbc$ is a finite $\bbc$-linear combination of
elements of $A(K_i)\otimes\bbc$, where the $K_i$ are finite Galois extensions of $K_0$.
Therefore, the span of the orbit $\Gal(K_0^{\sep}/K_0)v$ is contained in 
$A(L_0)\otimes\bbc$ for some finite Galois extension $L_0$ of $K_0$.
By the Mordell-Weil theorem for finitely generated fields \cite{Ner},
this is a finite dimensional vector space.  The image of $\sigma$ in $\Gal(L_0/K_0)$ is
finite, so the action of $\sigma$ on $A(L_0)\otimes\bbc$ is semisimple and all eigenvalues
are roots of unity.
\end{proof}

The implication (1) implies (2)  of Theorem~\ref{Spectrum} is trivial.  The statement (3) implies (4) 
is obvious since a closed subgroup of finite index in a topological group is also open.
For (2) implies (3) we have:

\begin{lem}
Let $A$ be an abelian variety over $K_0$ and $n$ a positive integer.
Let $\sigma\in\Gal(K_0^{\sep}/K_0)$.  If $\lambda$, an eigenvalue of
$\sigma$ acting on $A(K_0^{\sep})\otimes \bbc$, is a primitive $n$th root of unity, then
$$|\overline{\langle\sigma\rangle}/\overline{\langle\sigma^n\rangle}| = n.$$
\end{lem}

\begin{proof}
If $v$ is an eigenvector of $A(K_0^{\sep})\otimes\bbc$ with an eigenvalue $\lambda$ which is
a root of unity of order $n$, then $v\in A(L_0)\otimes\bbc$ for some finite Galois extension
$L_0/K_0$.  As $\sigma^k v = \lambda^k v$, the order $N$ of the image of $\sigma$ in 
$\Gal(L_0/K_0)$ is a multiple of the order of $\lambda$, i.e., a multiple of $n$.
Therefore, the kernel $\Delta$ of the composition 
$$\overline{\langle\sigma\rangle}\hookrightarrow 
\Gal(K_0^{\sep}/K_0)\twoheadrightarrow\Gal(L_0/K_0)$$
is of index $N$.  The composition
$$\overline{\langle\sigma\rangle} \twoheadrightarrow \overline{\langle\sigma\rangle}/\Delta
= \bbz/N\bbz \twoheadrightarrow \bbz/n\bbz$$
is therefore surjective, and its kernel $\overline{\langle\sigma^n\rangle}$ has index $n$.
\end{proof}

Finally, (4) implies (1) depends on applying Theorem~\ref{Rank} 
to abelian varieties obtained from $A$ by
restriction of scalars:

\begin{prop}
Let $A$ be an abelian variety over $K_0$ and $n$ a positive integer.
Let $\sigma\in\Gal(K_0^{\sep}/K_0)$.   If $\overline{\langle\sigma\rangle}$ has an open subgroup of index $n$, then every primitive $n$th root of unity is an eigenvalue of infinite multiplicity for $\sigma$
acting on $A(K_0^{\sep})\otimes\bbc$.
\end{prop}

\begin{proof}
Let $K\subset K_0^{\sep}$ denote the fixed field of $\sigma$ and $A_K$ denote the variety obtained
from $A$ by extending scalars to $K$.  Let $K_i\subset K_0^{\sep}$
denote the fixed field of the open subgroup of $\overline{\langle\sigma\rangle}$ of index $i$ if
such a field exists.  In particular, $K_d$ exists for all divisors $d$ of $n$.
Let $A_{d,K}$ denote the abelian variety over $K$ obtained from $A_{K_d}$ by restriction of scalars
from $K_d$ to $K$.  That is, $A_{d,K}$ represents the functor from $K$-schemes to 
abelian groups given by
$$A_{d,K}(S) = A(S\otimes_K K_d).$$
If $d\vert n$, there is a natural inclusion homomorphism
$i_{d,n}\colon A_{d,K}\to A_{n,K}$ which corresponds 
to the morphism $A(S\otimes_K K_d)\to A(S\otimes_K K_n)$ associated to the
inclusion $K_d\hookrightarrow K_n$.  Let
$$A_{n,K}^{\prim} := A_{n,K} \Bigm/ \sum_{\substack{d\vert n\\ d < n \\ }} i_{d,n}(A_{d,K}).$$
By Theorem~\ref{Rank}, $A_{n,K}^{\prim}(K)$ is of infinite rank.   This means
$$A(K_n)\otimes\bbq \Bigm/ \sum_{\substack{d\vert n\\ d < n \\ }} A(K_d)\otimes\bbq$$
is an infinite dimensional $\bbq$-vector space.  By construction, no proper subgroup of
$\bbz/n\bbz$ has any invariants on this space.  It follows that it is a direct sum of spaces
on which $\sigma$ acts with all eigenvalues primitive $n$th roots of unity.  As this is a rational representation,
all eigenvalues occur with infinite multiplicity.
\end{proof}

\end{document}